\documentclass[11pt, reqno]{amsart}

\usepackage{amssymb}

\evensidemargin\oddsidemargin

\begin{document}
\baselineskip=18pt
\setcounter{page}{1}
    
\newtheorem{Ass}{Assumption}
\newtheorem{Asso}{Assumption A\!\!}
\newtheorem{Assa}{Assumption B\!\!}
\newtheorem{Lemm}{Lemma}
\newtheorem{Rem}{Remark}
\newtheorem{Coro}{Corollary}
\newtheorem{Propo}{Proposition}

\renewcommand{\theAsso}{}
\renewcommand{\theAssa}{}

\def\a{\alpha}
\def\b{\beta}
\def\B{{\bf B}} 
\def\C{{\bf C}} 
\def\cG{{\mathcal{G}}} 
\def\cH{{\mathcal{H}}} 
\def\cI{{\mathcal{I}}} 
\def\cS{{\mathcal{S}}}
\def\UU{{\mathcal{U}}}
\def\ca{c_{\a}}
\def\ka{\kappa_{\a}}
\def\coa{c_{\a, 0}}
\def\cua{c_{\a, u}}
\def\cL{{\mathcal{L}}} 
\def\cM{{\mathcal{M}}} 
\def\Ea{E_\a}
\def\eps{{\varepsilon}} 
\def\esp{{\mathbb{E}}} 
\def\Ga{{\Gamma}} 
\def\GG{{\bf \Gamma}} 
\def\e{{\rm e}}
\def\ii{{\rm i}}
\def\L{{\bf L}}
\def\lbd{\lambda}
\def\lacc{\left\{}
\def\lcr{\left[}
\def\lpa{\left(}
\def\lva{\left|}
\def\M{{\bf M}}
\def\NN{{\mathbb{N}}} 
\def\pb{{\mathbb{P}}}
\def\QQ{{\mathbb{Q}}} 
\def\R{{\bf R}}
\def\rl{{\mathbb{R}}}
\def\racc{\right\}}
\def\rpa{\right)}
\def\rcr{\right]}
\def\rva{\right|}
\def\W{{\bf W}}
\def\X{{\bf X}}
\def\XX{{\mathcal X}}
\def\YY{{\mathcal Y}}
\def\Y{{\bf Y}}
\def\V{{\bf V}_\a}
\def\Un{{\bf 1}}
\def\S{{\bf S}}
\def\A{{\bf A}}
\def\G{{\bf G}}
\def\AA{{\mathcal A}}
\def\hAA{{\hat \AA}}
\def\hL{{\hat L}}
\def\T{{\bf T}}

\newcommand*\pFqskip{8mu}
\catcode`,\active
\newcommand*\pFq{\begingroup
        \catcode`\,\active
        \def ,{\mskip\pFqskip\relax}%
        \dopFq
}
\catcode`\,12
\def\dopFq#1#2#3#4#5{%
        {}_{#1}F_{#2}\biggl[\genfrac..{0pt}{}{#3}{#4};#5\biggr]%
        \endgroup
}

\def\claw{\stackrel{d}{\longrightarrow}}
\def\elaw{\stackrel{d}{=}}
\def\qed{\hfill$\square$}

\title{On a moment problem related to Bernstein functions}

\author[Thomas Simon]{Thomas Simon}

\address{Laboratoire Paul Painlev\'e, Universit\'e de Lille,  Cit\'e Scientifique, 59655 Villeneuve d'Ascq Cedex, France. {\em Email}: {\tt simon@math.univ-lille1.fr}}



\begin{abstract}
We give a simple proof of the moment-indeterminacy of the sequence $(n!)^t$ for $t > 2,$ using Lin's condition. Under a logarithmic self-decomposability assumption, the method conveys to power sequences defined as the rising factorials of a given Bernstein function, and to more general infinitely divisible moment sequences. We also provide a very short proof of the infinite divisibility of all the moment sequences recently investigated in \cite{LinC}, including Fuss-Catalan's.
\end{abstract}

\maketitle
\section{An example with the Gumbel distribution}

Let $\L$ be the standard exponential random variable and $\G =\log \L.$ For all $s > -1,$ one has
\begin{equation}
\label{Malm}
\esp [e^{s\G}]\; =\; \esp[\L^s]\; =\; \Ga(1+s)\; =\; \exp\lacc -\gamma s\, +\, \int^0_{-\infty} (e^{sx} - 1 - sx)\, \frac{dx}{\vert x\vert (e^{\vert x\vert} -1)}\racc
\end{equation}
where $\gamma$ is Euler's constant - see e.g. Formul\ae\, 1.7.2(19) and 1.9(1) in \cite{EMOT} for the third equality. This shows that $\G$ is infinitely divisible. Let $\{\G_t, \, t\ge 0\}$ be the real L\'evy process starting at zero such that $\G_1\elaw\G.$ For every $s >-1$ and $t\ge 0,$ one has
\begin{equation}
\label{Malm2}
\esp[e^{s\G_t}]\; =\; \Gamma(1+s)^t.
\end{equation}
The family of positive random variables $\L_{t} = e^{\G_t},\, t\ge 0,$ induces a multiplicative convolution semi-group, that is 
$$\L_{u}\L_{t}^{-1}\,\perp\;\L_{t}\qquad\mbox{and}\qquad\L_{u}\L_{t}^{-1}\,\elaw\, \L_{u-t}$$
if $t\le u.$ For every $t > 0,$ it is easy to see by (\ref{Malm2}) and Fourier inversion that the random variable $\L_t$ has a smooth density $f_t$ on $(0,\infty).$ It can also be shown, using (\ref{Malm2}) and Mellin inversion, that this density satisfies the integral equation
$$f_t(x)\; =\; \frac{1}{\Ga(t)}\; \int_x^\infty f_t(y)\, \lpa\log y - \log x\rpa^{t-1}\, dy,$$
but we shall not need this in the sequel. The positive entire moments of $\L_t$ are
$$m_n(t) \; =\; \esp[\L_t^n]\; =\; (n!)^t$$
and the following observation was made in \cite{Berg1}:
\begin{equation}
\label{MDMA}
\L_t\;\mbox{is moment-determinate}\;\Leftrightarrow\; t\le 2.
\end{equation}
See the recent survey \cite{Lin} for more details and references on the classical moment problem. Throughout, we will use the unusual but short notations MD for moment-determinate and MI for moment-indeterminate. 

As shown in \cite{Berg1}, the if part of (\ref{MDMA}) is an immediate application of Carleman's criterion: if $t\le 2$ one has, by Stirling's formula,
$$\sum_{n\ge 1} m_n(t)^{-\frac{1}{2n}}\;\ge\;\sum_{n\ge 1} m_n(2)^{-\frac{1}{2n}}\; =\; \sum_{n\ge 1} (n!)^{-\frac{1}{n}}\; =\; \infty.$$
In \cite{Berg1}, the proof of the only if part of (\ref{MDMA}) is however much more involved - see Section 2 therein. It amounts to checking the classical Krein condition for a certain associate distribution. Another proof was recently given in \cite{BergP}, using Krein's condition directly and an asymptotic analysis of $f_t$ at infinity. \\

We begin this note with an alternative and simple argument for the only if part of (\ref{MDMA}), which relies on Lin's condition - see Condition $L$ in Section 5 of \cite{Lin}. We use the HCM property of Bondesson and Thorin, which we will not introduce here in detail for the sake of concision. See Chapters 5 and 6 in \cite{Bond} for all definitions and notations.

\begin{Propo}
\label{HCMLt}
One has
$$f_t\;\mbox{{\em is HCM}}\; \Leftrightarrow\; t\ge 1.$$
\end{Propo}

\proof

We begin with the only if part. Suppose $t < 1$ and set $\M_t = \L_t^{\frac{1}{t}}.$ One has
$$\frac{\esp[\M_t^n]^{\frac{1}{n}}}{n}\; =\; \frac{\Ga(1+\frac{n}{t})^{\frac{t}{n}}}{n}\; \to\; \frac{1}{t\e}$$
by Stirling's formula. By Lemma 3.2 in \cite{CSY}, this implies
$$\log\pb[\L_t > x]\; =\; \log\pb[\M_t > x^{\frac{1}{t}}]\; \sim\; -t x^{\frac{1}{t}}.$$
Since $t <1,$ the upper tails of $\L_t$ are superexponentially small and $\L_t$ is not infinitely divisible, as is well-known. In particular, its density $f_t$ cannot be HCM.\\

We now prove the if part, which is clear for $t=1$ since $f_1(x) = e^{-x}.$ If $t > 1,$ introduce for every $q > 0$ the associate random variables
$$\L_{q,t}\; =\; \T\lpa \frac{q}{q+t-1}, \frac{1}{q+t-1}, \frac{t}{q+t-1}\rpa$$ 
with the notation of \cite{LS2}. It follows from (2.2) in \cite{LS2} and an immediate asymptotic analysis that for every $s\ge 0,$ 
$$\esp[\L_{q,t}^s]\; \to \; c_t^s \exp\lacc t\int^0_{-\infty} (e^{sx} - 1 - sx)\, \frac{dx}{\vert x\vert (e^{\vert x\vert} -1)}\racc$$
as $q\to\infty,$ for a constant $c_t$ to be determined. The normalization 
$$\esp[\L_{q,t}]\; =\; \esp[\L]\; =\; 1$$
for all $q > 0$ and $t> 1$ implies $c_t = e^{-t\gamma}$ and we can deduce from (\ref{Malm}) and (\ref{Malm2}) that
\begin{eqnarray}
\label{conv}
\L_{q,t} & \claw & \L_t
\end{eqnarray}
as $q\to\infty.$ Moreover, it follows from (2.4) and (2.7) in \cite{LS2} that
\begin{eqnarray*}
\L_{q,t} & \elaw &  \T\lpa \frac{q}{q+t-1}, \frac{1}{q+t-1}, \frac{1}{q+t-1}\rpa\,\times\, \T\lpa \frac{q+1}{q+t-1}, \frac{1}{q+t-1}, \frac{t-1}{q+t-1}\rpa\\
& \elaw & \lpa\frac{q+t-1}{q}\rpa\,\times\,\GG_{\frac{q}{q+t-1}}\,\times\, \T\lpa \frac{q+1}{q+t-1}, \frac{1}{q+t-1}, \frac{t-1}{q+t-1}\rpa
\end{eqnarray*}
where $\GG_\lambda$ stands for the standard Gamma random variable with parameter $\lambda > 0,$ and the products on the right-hand sides are independent. Applying Lemma 1 in \cite{BS3} shows now that $\L_{q,t}$ has a HCM density for every $q > 0, t > 1.$ By (\ref{conv}) and Theorem 5.1.3 in \cite{Bond}, this is also the case for $\L_t.$

\endproof

\begin{Coro}
The random variable $\L_t$ is {\em MI} for all $t>2.$ 
\end{Coro}

\proof
 
If $t > 2,$ combining Proposition \ref{HCMLt} and Property v) p.68 in \cite{Bond} shows that the function
$$x\;\mapsto\; -\frac{x f_t'(x)}{f_t(x)}$$
increases on $(0,\infty),$ which is Lin's condition. Moreover, we have
$$\frac{m_{n+1}(t)}{m_n(t)}\; =\; (n+1)^t.$$
By Theorem 5 in \cite{Lin}, $\L_t$ is MI.

\endproof

\begin{Rem}{\em (a) In the terminology of \cite{Wata}, the equivalence (\ref{MDMA}) means that the moment-determinacy of $\L_t$ is a time-dependent property for the L\'evy process $\{\G_t,\, t\ge 0\}.$ This temporal change comes from the multiplicative character of the semi-group associated to $\L_t.$ 

\medskip

(b) The HCM property is sensitive to power transformations for a given random variable - see Chapter 5 in \cite{Bond}. It is well-known that this is also the case for the moment problem - see the references in \cite{Lin}. Lin's condition is implied by the HM property of Section 6.4 in \cite{Bond}, which is less stringent than the HCM property and does not depend on powers. Another pleasant feature of Lin's condition is that it only deals with the behaviour of the density at infinity, in accordance with Krein's condition. This will be used in the next section.}
\end{Rem}

\section{A generalization}

If $\Phi$ is a Bernstein function, it was proved in \cite{BY1} that the sequence
$$\Phi(1)\,\times\,\cdots\,\times\,\Phi(n)$$
is a determinate moment sequence on $\rl^+.$ The corresponding positive random variable $\R$ is a multiplicative factor of $\L$ which is called in \cite{HY} the Remainder. In Theorem 1.8 of \cite{Berg1}, it is shown that $\S = \log \R$ is infinitely divisible - see also Theorem 3.1 in \cite{HY} for a different proof. This leads to a real L\'evy process $\{\S_t, \, t\ge 0\}$ characterized by $\S_1 =\S,$ and to the family of positive random variables $\R_t  = e^{\S_t},\, t\ge 0.$ The latter induces a multiplicative convolution semi-group exactly as above. The positive entire moments of $\R_t$ are
$$\mu_n(t) \; =\; \esp[\R_t^n]\; =\; (\Phi(1)\,\times\,\cdots\,\times\,\Phi(n))^t,\qquad n\ge 0.$$
The unit drift case $\Phi(x) = x$ yields $\R_t =\L_t,$ which is for $t\in (0,1)$ the Remainder associated to a stable subordinator of parameter $t$ - see Example 4.1 in \cite{HY}. In view of the previous section, it is natural to ask for the MD or MI character of $\R_t.$ This problem was recently adressed in \cite{PV}. Introduce the following parameter
$$\ell\; =\; \limsup_{x\to\infty}\lpa \frac{\Psi(x)}{x \,\log x}\rpa\quad\mbox{with}\quad\Psi(x)\; =\; \int_0^x \log \Phi (t)\, dt.$$
Observe from the concavity of $\Phi$ that necessarily one has $\ell\in [0,1].$ We begin with an easy and general result.

\begin{Propo}
\label{EaG}
The random variable $\R_t$ is {\em MD} for all $t<2/\ell.$ 
\end{Propo}

\proof

Since $\log\Phi$ increases, one has
$$\ell\; =\; \limsup_{n\to\infty}\lpa \frac{\log(\Phi(1)\times\cdots\times\Phi(n))}{n \,\log n}\rpa.$$
Therefore,
$$\limsup_{n\to\infty} \lpa\frac{\log\mu_n(t)}{n \,\log n}\rpa\; < \;2$$

\smallskip

\noindent
whenever $t < 2/\ell,$ whence the conclusion by Carleman's criterion.

\endproof

In particular, we see that $\R_t$ is MD for every $t > 0$ if $\ell =0,$ which is the case in Examples 4.2, 4.3 and 4.4 of \cite{HY}. In the following, we implicitly assume $\ell > 0.$ 

\medskip

It is natural to ask for the MI character of $\R_t$ when $t > 2/\ell.$ A first difficulty is that in order to check either Krein's or Lin's condition, the absolute continuity of the law of $\R_t$ is required. The latter is equivalent to that of $\S_t$ but the problem of absolute continuity for marginals of a real L\'evy process is hard in general, subject to temporal changes - see \cite{Wata}. We will consider the following 

\begin{Ass} The random variable $\S$ is self-decomposable.
\end{Ass}

Under this assumption, it is well-known that $\S_t$ is absolutely continuous for every $t > 0.$ The self-decomposability of $\S$ can be characterized by the Bernstein measure $\kappa(dx)$ of the completely monotonic function $\Phi'/\Phi.$ It is shown in Proposition 3.5 of \cite{HY} that $\S$ is self-decomposable if and only if $\kappa$ is absolutely continuous with a density $\kappa$ such that $x\mapsto (e^x-1)^{-1} \kappa(x)$ is non-increasing on $(0,\infty).$ This is true for Examples 4.1, 4.5 and 4.6 of \cite{HY}. This is also true if $\Phi$ is a complete Bernstein function, because it admits the representation
$$\Phi(x)\; =\; \Phi(1)\,\exp\lacc \int_0^\infty \lpa \frac{1}{1+t} - \frac{1}{x+t}\rpa \eta(t)\, dt\racc$$
for some measurable function $\eta$ taking its values in $[0,1]$ - see e.g. Theorem 6.10 in \cite{SSV}. A simple computation shows then that the measure $\kappa$ has density
$$\kappa(x)\; =\; x\int_0^\infty e^{-xt}\, \eta(t)\, dt,$$
and the function $(e^x -1)^{-1} \kappa(x)$ is hence non-increasing on $(0,\infty).$ Let us also mention that by Corollary 1.11 of \cite{AJR}, the random variable $\S$ is self-decomposable when the upper tail of the L\'evy measure of $\Phi$ is log-convex, which is less stringent than the complete Bernstein character of $\Phi.$ See \cite{AJR} and the references therein for further aspects of the measure $\kappa.$ Introducing the further parameter
$${\bar \ell}\; =\; \liminf_{x\to\infty}\lpa \frac{\Psi(x)}{x \,\log x}\rpa,$$
we have the following counterpart to Proposition \ref{EaG}.
 
\begin{Propo}
\label{EaF}
Under Assumption {\em 1}, the random variable $\R_t$ is {\em MI} for every $t>2/{\bar \ell}.$ 
\end{Propo}

\proof
If $t > 2/{\bar \ell}$ we obtain, reasoning as in Proposition \ref{EaG},
$$\liminf_{n\to\infty} \lpa\frac{\log\mu_n(t)}{n \,\log n}\rpa\; > \;2.$$
Hence, by Theorem 7 in \cite{Lin}, we just need to check Lin's condition on $\R_t.$ Recall now that the L\'evy measure of the self-decomposable random variable $-\S_t$ has, by Proposition 3.2 in \cite{HY}, density 
$$t\,x^{-1}(e^x-1)^{-1} \kappa(x)\,\Un_{(0,\infty)}(x).$$
Moreover, it follows from Lemma 1.3 in \cite{AJR} that
$$\log\Phi(x)\; =\; \log\Phi(1)\; +\; \int_0^\infty (e^{-t} - e^{-xt})\, \frac{\kappa(t)}{t}\, dt.$$ 
Therefore, one has necessarily 
$$\int_0^1 (e^x-1)^{-1} \kappa(x)\, dx \; =\; \infty$$
since otherwise $\Phi$ would be bounded, which is clearly excluded by the positivity of $\ell.$ All of this shows that the law of $-\S_t$ is of the type ${\rm I}_7$ in \cite{SY} and, by Theorem 1.3 (xii) therein, that the density of $\S_t$ is log-concave on some interval $[a,\infty).$ Changing the variable, the latter is equivalent to Lin's condition on $\R_t.$

\endproof

\begin{Rem}
\label{EaF+}
{\em (a) The monotone density theorem applied to the concave function $\log\Phi$ shows that
$$\ell\, =\, {\bar \ell}\;\Leftrightarrow\; \log \Phi(x)\, \sim\, \ell\,\log x.$$
The latter holds true when $\Phi$ is regularly varying at infinity, the index being then necessarily $\ell.$ Observe also that this is a weaker condition than regular variation, given at the logarithmic level.

\medskip

(b) If $\ell ={\bar \ell},$ one may wonder if $\R_{2/\ell}$ is MD or MI. The behaviour at the threshold is usually a question where trouble begins in the literature on moment problems, in the absence of universal criteria. Under Assumption 1, a consequence of the above proof and Theorem 3 in \cite{Pak} is
$$\R_t \;\mbox{is MD}\quad\Leftrightarrow\quad \sum_{n\ge 1}\, \mu_t(n)^{-\frac{1}{2n}}\; =\; \infty.$$
In particular, one has
\begin{itemize}
\item $\R_{2/\ell}$ is MI if
$$\liminf_{x\to\infty}\lpa \frac{\Phi(x)}{x^\ell \,(\log x)^c}\rpa\; =\; \infty\qquad\mbox{for some $c >1.$}$$
\item $\R_{2/\ell}$ is MD if
$$\limsup_{x\to\infty}\lpa \frac{\Phi(x)}{x^\ell \,(\log x)^c}\rpa\; <\; \infty\qquad\mbox{for some $c <1.$}$$
\end{itemize}

\smallskip

\noindent
See Remark 3 below for another, Krein type, criterion.

\medskip

(c) In the unit drift case $\Phi(x)=x$, Assumption 1 is fulfilled and Proposition \ref{EaF} gives another quick proof of the MI character of the moment sequence $(n!)^t$ for $t > 2.$ Our previous argument in Section 1 is more involved, but it is also more informative.}

\end{Rem}

\medskip

The above proof enhanced Lin's condition. We now show that it is possible to derive an analogous result with the help of the classical Krein's condition, under two different assumptions. The first one is weaker than Assumption 1.

\begin{Asso}
There exists $\varepsilon > 0$ such that for every $t \in (2/\ell, 2/\ell +\varepsilon),$ the random variable $\R_t$ has a density which is ultimately monotone.
\end{Asso}

By Yamazato's theorem on the unimodality of self-decomposable laws and a change of variable, Assumption A is implied by Assumption 1. The study of monotonicity properties of ID densities on $\rl$ can be a delicate problem, leading to pathological situations. For example, an ID density may have an infinite number of modes - see again \cite{Wata} for more on this topic. We were not able to exhibit any absolutely continuous Remainder not fulfilling Assumption A, but we believe there should be some. 

\begin{Assa} 
There exists $c \in (0,\infty)$ such that 
$$\lim_{n\to\infty} \frac{\lpa \Phi(1)\,\times\,\cdots\,\times\,\Phi(n)\rpa^{\frac{1}{n\ell}}}{n}\; =\; c.$$
\end{Assa}

\medskip

By monotonicity, this second assumption implies $\ell = {\bar \ell}.$ It is also a stronger condition, given at the natural level. Observe finally that Assumption B is fulfilled when $\Phi$ is regularly varying at infinity, the index being then necessarily $\ell.$
  
\begin{Propo}
\label{True}

Under Assumptions {\em A} and {\em B}, the random variable $\R_t$ is {\em MI} for every $t>2/\ell.$ 
\end{Propo}

\proof

As is well-known, it is enough to consider the case $t\in(2/\ell, 2/\ell +\varepsilon).$ Fix such a $t$ and set 
$$\nu_n(t)\; =\; \esp[\R_t^{\frac{n}{\ell t}}]$$ 
for all $n\ge 1.$ Since $p\mapsto \esp[X^p]^{\frac{1}{p}}$ is non-increasing on $(0,\infty)$ for any positive random variable $X,$ one has the bounds 
$$\frac{[n(\ell t)^{-1}]}{n}\,\times\,\frac{\esp[\R^{[n(\ell t)^{-1}]}]^{\frac{1}{\ell[n(\ell t)^{-1}]}}}{[n(\ell t)^{-1}]} \;\le\; \frac{\nu_n(t)^{\frac{1}{n}}}{n}\; \le\; \frac{[n(\ell t)^{-1}] +1}{n}\,\times\,\frac{\esp[\R^{[n(\ell t)^{-1}]+1}]^{\frac{1}{\ell[n(\ell t)^{-1}]+1}}}{[n(\ell t)^{-1}]+1}$$
and Assumption B implies
$$\lim_{n\to\infty} \lpa  \frac{\nu_n(t)^{\frac{1}{n}}}{n}\rpa\; =\; \frac{c}{\ell t}\cdot$$
Reasoning as in Proposition 1 leads then to the estimate
$$\log\pb[\R_t > x]\; \sim\; - c_t\, x^{\frac{1}{\ell t}}$$
with $c_t = \ell t(\e c) ^{-1}\in (0,\infty).$ By Assumption A, the random variable $\R_t$ has a density $f_t$ which is non-increasing at infinity, and the above estimate implies easily 
$$x^2 f_t(x)\; \ge \; e^{-(c_t/2) x^{\frac{1}{\ell t}}}$$
for $x$ large enough. Since $t > 2/\ell,$ the relaxed Krein's condition given in Theorem 4 of \cite{Lin} is in force, and $\R_t$ is MI.

\endproof

\section{Infinitely divisible moment sequences}

If $\{\mu_n, \, n\ge 1\}$ is the entire moment sequence of a positive random variable $\X,$ we say that this sequence is infinitely divisible (ID for short) if $\{\mu_n^t, \, n\ge 1\}$ is an entire moment sequence for every $t > 0.$ It is clear from the considerations in Section 1 that this property is equivalent to the infinite divisibility of $\log \X$ as a random variable. In particular, every moment sequence
$$\mu_n\; =\; \Phi(1)\;\times\;\cdots\;\times\;\Phi(n)$$
with $\Phi$ a Bernstein function, is ID - see the aforementioned Theorem 1.8 in \cite{Berg1} and Theorem 3.1 in \cite{HY}. For short, we will say that such an ID moment sequence is Bernstein. Observe that ID moment sequences need not be Bernstein, as shows the example $\mu_n = (n!)^2.$ 

An entire moment sequence $\{\mu_n, \, n\ge 1\}$ which is both ID and MD gives rise to two sets of positive random variables indexed by time, via the associate random variable $\X.$ The first one is the multiplicative family $\{\X_t, \, t > 0\}$ defined as in Section 1, and the second one is the family of positive power transformations $\{\X^t, \, t > 0\}.$ A conjecture formulated in \cite{LinC} is that $\X_t$ is MD if and only if $\X^t$ is MD, for every $t >0.$ The following proposition gives a partial answer. Recall that the L\'evy measure of a real self-decomposable random variable has density $k(x)/\vert x\vert$ on $\rl^*,$ where $k$ is a function non-decreasing on $(-\infty,0)$ and non-increasing on $(0,\infty),$ which is called the spectral function. Henceforth, we implicitly exclude the case where $\log\X$ has a non-trivial Gaussian component, since then both $\X_t$ and $\X^t$ have the log-normal distribution as multiplicative factor for every $t >0$ and are hence MI.  

\begin{Propo}
\label{LinA}
Assume that the random variable $\log \X$ is self-decomposable and that its spectral function is not integrable at $0-$. Then, for every $t>0,$ one has
$$\X_t\;\mbox{is {\em MD}}\quad\Leftrightarrow\quad \X^t\;\mbox{is {\em MD}}.$$
\end{Propo}

\proof

The assumption means that $-\log\X$ is of the type ${\rm I}_7$ in \cite{SY}, which implies as in the proof of Proposition 3 that Lin's condition is satisfied by both $\X_t$ and $\X^t.$ Applying Theorem 3 in \cite{Pak}, we get
$$\X_t\;\mbox{is MD}\quad\Leftrightarrow\quad\sum_{n\ge 1} \mu_n^{-\frac{t}{2n}}\; =\; \infty\quad\Leftrightarrow\quad\sum_{n\ge 1} \mu_{[nt]}^{-\frac{t}{2[nt]}}\; =\; \infty\quad\Leftrightarrow\quad \X^t\;\mbox{is MD},$$
the second equivalence being an easy consequence of the non-increasing character of $n\mapsto \mu_n^{-\frac{t}{2n}},$ whereas the third equivalence is obtained as for the bounds in the proof of Proposition 4.
 
\endproof

\begin{Rem}{\em For a Remainder $\R$ satisfying Assumption 1, Proposition \ref{LinA} combined with Theorems 4 and 10 in \cite{Lin}, and a change of variable, shows that at the threshold $t = 2/\ell,$ one has the Krein type criterion
$$\R_{2/\ell}\;\mbox{is MD}\quad\Leftrightarrow\quad \R^{2/\ell}\;\mbox{is MD}\quad\Leftrightarrow\quad \int_0^\infty \frac{-\log f(x^\ell)}{1+x^2}\, dx\; =\; \infty$$
where $f$ is the density of $\R.$ This is useful since the asymptotic analysis at infinity of the density of $\R_t$ for $t\neq 1$ might be more involved than that of $f$ - see \cite{BergP} for the case $\R =\L.$}
\end{Rem}

In the recent paper \cite{LinC}, the infinite divisibility of several classical moment sequences was obtained for the first time, using essentially Theorem 1.8 in \cite{Berg1}. We can show this property for two larger families of moment sequences, and in a very simple way. We will also investigate some other interesting moment sequences.

\subsection{Gamma moment sequences of order 2} We consider the moment sequence 
$$\mu_n\; =\; \frac{\Ga(a+sn)\Ga(a+b)}{\Ga(a)\Ga(a+b +sn)}$$
for $a,b,s > 0.$ The associate random variable is the power transformation $\B_{a,b}^s$ of the standard Beta random variable $\B_{a,b}$ with density
$$\frac{\Ga(a+b)}{\Ga(a)\Ga(b)}\, x^{a-1} (1-x)^{b-1}\, \Un_{(0,1)} (x).$$ 
This sequence is ID by the well-known fact that $\log\B_{a,b}^s = s\log\B_{a,b}$ is an ID random variable. The latter amounts to the standard Malmst\'en type formula
$$\esp[e^{\lbd\log\B_{a,b}^s}]\; =\; \frac{\Ga(a+s\lbd)\Ga(a+b)}{\Ga(a)\Ga(a+b +s\lbd)}\; =\; \exp\lacc-\int_0^\infty (1-e^{-\lbd x}) \lpa\frac{e^{-ax}(1-e^{-bx})}{x (1-e^{-x})}\rpa dx\racc.$$
Clearly, one has $\mu_n^t\to 0$ as $n\to \infty$ so that $\{\mu_n^t\}$ is MD for every $t > 0,$ by Carleman's criterion. The same is true for $\{\mu_{nt}\}$ since the associate random variable $\B_{a,b}^{st}$ has bounded support.

\medskip

Taking $a=1/2, b=3/2$ and $s=1,$ we have
$$\mu_n\; =\; \frac{1}{4^n(n+1)}\,\binom{2n}{n}\; =\; 4^{-n}\,C_n$$
which is the Catalan number sequence up to a multiplicative constant. Hence, the previous discussion encompasses Theorem 1 in \cite{LinC}. Moreover, letting $b\to\infty,$ we get
$$b^{sn}\mu_n\; \to\; \frac{\Ga(a+sn)}{\Ga(a)}$$
for every $a,s > 0,$ which is the Gamma sequence of order 1 recently studied in \cite{Berg2}. Since the ID property of a moment sequence is preserved under pointwise limit, the latter sequence is also ID. The associate random variable is the power transformation $\GG_a^s$ of the standard Gamma random variable $\GG_a$ with density
$$\frac{1}{\Ga(a)}\; x^{a-1} e^{-x}\, \Un_{(0,\infty)} (x),$$
which satisfies Lin's condition. Hence, Carleman's criterion applies in both directions and Stirling's formula shows that the sequence
$$\lpa\frac{\Ga(a+sn)}{\Ga(a)}\rpa^t$$ 
is MD if and only if $st\le 2.$ Putting everything together, we have got a very simple proof of Conjecture 2 in \cite{LinC}, which is also Theorem 1.1 in \cite{Berg2}. 
 
\begin{Rem}
\label{Bern}
{\em (a) A natural and more involved question, which is connected to the approach of \cite{LinC}, is whether the sequence $\{\mu_n\}$ is Bernstein, with our above notation. This question was actually already adressed in \cite{BS1}, for other purposes. It follows easily from the hypergeometric transformations carried out in Section 2.3 of \cite{BS1} that
$$\{\mu_n\}\;\;\mbox{is Bernstein}\;\Leftrightarrow\; \inf\{b,s\}\,\le\, 1\quad\mbox{and}\quad a\,\ge\, s.$$
The corresponding Bernstein function is
$$\Phi(\lbd)\; =\; \frac{\Ga(a)\Ga(a+b-s)}{\Ga(a+b)\Ga(a-s)}\; +\; \int_0^\infty (1-e^{-\lbd x})\,\rho(x)\, dx$$ 
with 
$$\rho(x)\; =\; be^{-as^{-1}x}\;\pFq{2}{1}{1+s,1-b}{2}{1-e^{-s^{-1}x}}\; =\; be^{(1-(a+b)s^{-1})x}\;\pFq{2}{1}{1+b,1-s}{2}{1-e^{-s^{-1}x}},$$
a non-negative integrable function which simplifies into $\rho(x) = be^{(1- (a+b)s^{-1})x}$ for $s=1$ and into $\rho(x) = be^{-as^{-1}x}$ for $b=1.$ We refer to Section 2.3 in \cite{BS1} for some other interesting aspects of the function $\rho,$ connected to the zeroes of the classical hypergeometric series. In particular, it can be shown that $\Phi$ is a complete Bernstein function if and only if $b=1$ or $s=1,$ and that it belongs to the Jurek class if and only if $2a+b +s +bs \ge 1.$

\medskip

(b) For the Catalan moment sequence, one has $C_n = \Phi(1)\times\cdots\times\Phi(n),$ where
$$\Phi(\lbd) \; =\; 2\lpa 2\, -\frac{3}{1+\lbd}\rpa$$
is not Bernstein in the strict sence, since it takes negative values on $(0,1/2).$ On the other hand, the function ${\tilde \Phi}(\lbd) = \Phi(\lbd +1/2)$ is Bernstein and the factorization $C_n = {\tilde \Phi}(1/2)\times\cdots\times{\tilde \Phi}(n-1/2)$ was used in \cite{LinC} together with a previous Lemma of \cite{BergP} to show the ID character of $\{C_n\}$.

\medskip

(c) Similarly as above, one can show that
$$\lacc \frac{\Ga(a+sn)}{\Ga(a)}\racc\;\;\mbox{is Bernstein}\;\Leftrightarrow\; \inf\{1,a\}\,\ge\, s.$$
The corresponding Bernstein functions are
$$\Phi(\lbd)\; =\; \frac{\Ga(a)}{\Ga(a-s)}\; +\; \frac{1}{\Ga(1-s)}\int_0^\infty (1-e^{-\lbd x})\lpa \frac{e^{-as^{-1}x}}{(1-e^{-s^{-1}x})^{1+s}}\rpa dx$$ 
for $s < 1$ and $\Phi(\lbd) = a-1 +\lbd$ for $s=1.$}
\end{Rem}

\subsection{Binomial and Raney moment sequences} We consider the sequence 
$$\mu_n\; =\; \binom{pn +r}{n}$$
which is known - see \cite{MP} - to be a moment sequence on $\rl^+$ if and only if $p\ge 1$ and $r\in[-1,p-1].$ The associate random variable is here more complicated than above. It follows from Theorem 3.1 in \cite{MP} that for $p\ge 1$ rational and $r \in(-1,p-1],$ it is a renormalized finite product of independent Beta random variables, so that its logarithm is infinitely divisible. A density argument shows then immediately that the sequence $\{\mu_n\}$ is ID for all $p\ge 1$ and $r\in[-1,p-1].$ Carleman's criterion implies, as above, that the power sequences $\{\mu_n^t\}$ and $\{\mu_{nt}\}$ are MD for every $t>0$.

\medskip

Taking now $r=0$ and $p\ge 2$ an integer, we obtain a quick proof of Theorems 2 and 6 in \cite{LinC}, and also of Theorem 2', 3 and 5 therein by the factorization argument given in Lemma 2 of \cite{LinC}. For instance, the Fuss-Catalan sequence of order $k$
$$C_{k,n}\; =\; \lpa\frac{1}{1 +kn}\rpa\,\times\,\binom{(k+1)n}{n}$$   
is the product of two ID moment sequences taking $p = k+1, r= 0$ and, in the previous paragraph, $a=b=1, s=k.$
 
\medskip

In the same vein, it is interesting to mention that the above discussion also implies the ID character of the Raney sequence
$$\mu_n\; =\; \frac{r}{np+r}\,\binom{pn +r}{n},$$  
which is a moment sequence on $\rl^+$ if and only if $p\ge 1$ and $r\in[0,p]$ - see \cite{MP,MPZ}. Indeed, in the non-trivial case $r\neq 0,$ we have the factorization
$$\mu_n\; =\; \lpa\frac{1}{1 +r^{-1}(p-1)n}\rpa\,\times\,\binom{pn +r-1}{n}$$   
and we can again apply Lemma 2 in \cite{LinC}, the first factor being ID by the case $a=b=1$ and $s = r^{-1}(p-1)$ of the previous paragraph. Taking $p=2$ and $r=1,$ we recover the ID character of the Catalan sequence. 

\begin{Rem}
\label{BernR}
{\em Characterizing the Bernstein property of the binomial and the Raney moment sequences is an open problem, which is apparently not easy. Indeed, in both situations the ratio
$$\frac{\mu_n}{\mu_{n-1}}$$
involves six Gamma functions in general. Trying an hypergeometric summation argument as in Section 2.3 in \cite{BS1} should lead to generalized hypergeometric series, and it is well-known that exact summation formul\ae\, are rather rare in this broader context.}
\end{Rem}

\subsection{Gamma moment sequences of higher order} Consider the general sequence
$$\mu_n\; =\; \prod_{i=1}^p \frac{\Ga(a_i+A_i\, n)}{\Ga(a_i)}\;\times\, \prod_{j=1}^q \frac{\Ga(b_j)}{\Ga(b_j +B_j\, n)}$$
with all parameters positive. Characterizing the positive definiteness in $\rl^+$ of this sequence, that is whether it is the moment sequence of a positive random variable, as was done in \cite{MP, MPZ} for binomial and Raney sequences, seems to be a difficult task which has not been undertaken as yet. On the other hand, the ID moment sequence character of $\{\mu_n\}$ in the case of a compact support can be characterized from the recent results in \cite{KP16}. To be more precise, it follows from Lemma 1 and Theorem 4 in \cite{KP16} that $\{\mu_n\}$ is the moment sequence of a positive random variable $\X$ with compact support and that this sequence is ID, if and only if 
$$\sum_{i=1}^p A_i\; =\; \sum_{j=1}^q B_j\qquad\mbox{and}\qquad \sum_{i=1}^p \frac{e^{-a_iA_i^{-1} x}}{1- e^{-A_i^{-1} x}}\; -\; \sum_{j=1}^q \frac{e^{-a_jB_j^{-1} x}}{1- e^{-B_j^{-1} x}}\; \ge\; 0\quad\forall\, x\, > \, 0.$$
The support of $\X$ is then the interval $[0,\rho]$ with 
$$\rho\; =\; \prod_{i=1}^p A_i^{A_i}\,\times\, \prod_{j=1}^q B_j^{-B_j},$$ 
see (16) of \cite{KP16}. However, as mentioned in the introduction to \cite{KP16}, the non-negativity condition is not easy to check directly on the parameters, even in the special cases of Raney and binomial sequences. In these two cases, the above density argument via the Beta distribution is much quicker.

\begin{Rem}
\label{Ostrov}
{\em Other ID moment sequences with compact support can be built on two sequences of positive numbers via the multiple Gamma function. The associate random variables have the so-called Barnes Beta distribution - see \cite{Ost}, especially Theorem 2.4 therein, for more detail.}
\end{Rem}

\subsection{Other moment sequences of Gamma type} In this last paragraph we come back to the moment sequence $(n!)^t$ of Section 1. The third equality in (\ref{Malm}) yields the exponential representation
$$\frac{\Ga(1+s)^t}{\Ga(1+st)}\; =\; \exp\lacc \int^0_{-\infty} (e^{sx} - 1 - sx)\,\lpa \frac{t}{e^{\vert x\vert} -1} - \frac{1}{e^{\vert x\vert t^{-1}} -1}\rpa \frac{dx}{\vert x\vert}\racc$$
for every $s, t > 0.$ Besides, an elementary analysis shows that the function
$z\mapsto t(z-1) - (z^t-1)$ is positive on $(1,\infty)$ for $t\in(0,1)$ and negative on $(1,\infty)$ for $t > 1.$ By Mellin inversion, this shows the identities in law
$$\L_t\;\elaw\; \L^t\; \times\;\M_t\quad\mbox{for $t\in(0,1)$}\qquad\mbox{and}\qquad\L^t\;\elaw\; \L_t\; \times\;\M_t\quad\mbox{for $t>1,$}$$
where $\M_t$ is a positive random variable with fractional moments
$$\esp[\M_t^s]\; =\;\exp\lacc \int^0_{-\infty} (e^{sx} - 1 - sx)\,\lva \frac{t}{e^{\vert x\vert} -1} - \frac{1}{e^{\vert x\vert t^{-1}} -1}\rva \frac{dx}{\vert x\vert}\racc, \qquad s\, >\, -\inf\{1,t^{-1}\}.$$
Stirling's formula implies that $\M_t$ has a compact support which is $[0,t^{-t}]$ for $t\in(0, 1)$ and $[0,t^t]$ for $t > 1.$ The positive entire moments of $\M_t$ are
$$\mu_n\;=\; \frac{(n!)^t}{\Ga(1+nt)}\quad\mbox{for $t<1$}\qquad\mbox{and}\qquad\mu_n\; =\; \frac{\Ga(1+nt)}{(n!)^t} \quad\mbox{for $t>1.$}$$
Since $\log \M_t$ is infinitely divisible, these moment sequences are ID. Moreover, by compactness of the support, the sequences $\{\mu_{ns}\}$ and $\{\mu_n^s\}$ are MD for every $s >0.$ For $t > 1,$ the sequence $\{\mu_n\}$ is not Bernstein since the corresponding Remainder would then have Laplace exponent
$$\frac{\Ga(1+\lbd t)}{\lbd^t\Ga(1-t + \lbd t)}$$
and this function takes negative values on $(0,\infty)$. For $t \in(0, 1),$ the Bernstein character of $\{\mu_n\}$ amounts to that of the function
$$\Phi(\lbd)\; = \;\frac{\Ga(1-t+\lbd)}{\lbd^{1-t}\,\Ga(\lbd)}$$
but I was not able to give an answer to this interesting question. Combining Formul\ae\, 1.7.2(22) and 1.9(1) in \cite{EMOT} yields the exponential representation
$$\Phi(\lbd)\; =\; \exp\lacc \int_0^{\infty} e^{-\lbd x}\lpa \frac{(1-t)(1-e^{-x}) +e^{-(1-t)x} -1}{x(1-e^{-x})}\rpa dx \racc$$
and it can be shown, using Theorem 6.10 in \cite{SSV}, that the function on the right-hand side is not a complete Bernstein function. This representation also shows that $1/\Phi(\lbd)$ is logarithmically completely monotone, which is necessary but not sufficient for $\Phi(\lbd)$ to be Bernstein - see Proposition 5.17 in \cite{SSV} and the remark thereafter.

\section{A further example}

We conclude this paper with another example of strict dichotomy between moment-determinacy and moment-indeterminacy, in the spirit of Section 1. The framework is that of the r-gstable$(a,m)$ laws recently studied in \cite{Pakes, JSW} - see also the references therein. These laws are well-defined if and only if $0<a<m$ and they form a generalization of the inverse positive stable laws which correspond to the case $m=1.$ We refer to \cite{Pakes,JSW} for more details. The entire moments of the corresponding random variable $\Y_{a, m}$ are given by
$$\mu_n\; =\; \esp[\Y_{a,m}^n]\; =\;a^{\frac{(m-a)n}{a}}\times\, \frac{G(m+n,a) G(a,a)}{G(a+n,a)G(m,a)}$$

\smallskip

\noindent
where $G$ is the double Gamma function - see (11) in \cite{JSW}. It follows from the main Theorem in \cite{JSW} that $\log \Y_{a,m}$ is infinitely divisible and the moment sequence $\{\mu_n\}$ is hence always ID. It is also easy to see from the concatenation formula $G(z+1,\tau) = \Ga(z\tau^{-1})G(z,\tau)$ that
$$\{\mu_n\}\;\;\mbox{is Bernstein}\;\Leftrightarrow\; 1\,\le\, a\, <\, m\,\le 3a-1.$$
The corresponding Bernstein functions are
$$\Phi(\lbd)\; =\; a^{\frac{m-a}{a}}\lpa \frac{\Ga(\frac{m}{a}-1)}{\Ga(1-\frac{1}{a})}\; +\; \frac{m+1-2a}{\Ga(3-\frac{m+1}{a})}\int_0^\infty (1-e^{-\lbd x})\lpa \frac{e^{-(m-a)x}}{(1-e^{-ax})^{\frac{m+1}{a}-1}}\rpa dx\rpa$$ 
for $m < 3a -1$ and $\Phi(\lbd) = a^{1-\frac{1}{a}}(a-1 +\lbd)$ for $m =3a-1.$

\begin{Propo}
\label{MomY}
The random variable $\Y_{a,m}$ is {\em MD} if and only if $m\le 3a.$ 
\end{Propo} 

\proof An easy consequence of (11) and (13) in \cite{JSW} is
$$\esp[\Y_{a,m}^n]^{\frac{1}{n}}\; \sim\; \lpa\frac{n}{a\e}\rpa^{\frac{m-a}{a}},$$
so that 
$$\frac{\log\esp[\Y_{a,m}^n]}{2n\log n}\;\to\; \frac{m-a}{2a}\cdot$$
Moreover, we know by Corollary (b) in \cite{JSW} that $\Y_{a,m}$ has a HCM density for $m\ge 2a.$ A combination of Carleman's criterion and Theorem 7 in \cite{Lin} shows that $\Y_{a,m}$ is MD if and only if 
$$\frac{m-a}{2a}\,\le \, 1\;\Leftrightarrow\; m\le 3a.$$ 
\endproof

This result was obtained in Theorem 8.2 of \cite{Pakes} for $m$ integer, the only if part being there a consequence of Krein's condition and the subexponential tail behaviour at infinity of the density of $\sqrt{Y_{a,m}},$ which is obtained by means of a certain class of special functions. It is shown in the Proposition of \cite{JSW} that the latter subexponentiality property holds true for every $m >3a$ non necessarily an integer, so that one can conclude as in \cite{Pakes}. Overall, this proof is more involved than the above HCM argument.

\end{document}